\documentclass[12pt]{article}
\usepackage[leqno]{amsmath}
\usepackage{amssymb}
\usepackage{latexsym}

\newtheorem{theorem}{Theorem}[section]
\newtheorem{remark}[theorem]{Remark}

\newtheorem{proposition}[theorem]{Proposition}
\newtheorem{lemma}[theorem]{Lemma}
\newtheorem{corollary}[theorem]{Corollary}

\newtheorem{definition}[theorem]{Definition}

\numberwithin{equation}{section}

\def\Z{{ \mathbb Z}}

\def\End{{\rm End\,}}
\def\Hom{{\rm Hom\,}}

\def\cf{{\rm cf\,}}

\def\dom{{\rm dom\,}}
\def\Im{{\rm Im\,}}

\def\ra{\rightarrow}
\def\arr{\longrightarrow}

\begin{document}
\def\Ines{{\rm Ines\,}}
\def\Fin{{\rm Fin\,}}
\def\Small{{\rm Small\,}}

\title{ENDOMORPHISM RINGS OF MODULES WHOSE CARDINALITY IS COFINAL TO
OMEGA
\thanks {This work is supported by the project No. G-0294-081.06/93 of the 
German-Israeli Foundation for Scientific Research \& Development. 
}}

\author{R\"udiger G\"obel\\Fachbereich 6,\\ Mathematik und 
Informatik\\Universit\"at Essen, \\ 45117 Essen, Germany 
\\e-mail: {\small R.Goebel@Uni-Essen.De}
\and
Saharon Shelah
\thanks{Manuscript No. GbSh547}
 \\ Institute of Mathematics \\ Hebrew University\\
91904 Jerusalem \\ Israel 
\\e-mail: {\small Shelah@math.huji.ac.il}}
\date{}

\maketitle



\section{Introduction}

We want to consider torsion-free $R$-modules over a ring $R$. In Section 3
the ring $R$ will be a principal ideal domain and in Section 4 we allow
more general commutative rings $R$. However generally we assume that
$R$ has a distinguished countable, multiplicatively closed 
subset $S$ of non-zero
divisors. We also may assume that $1 \in S$ and say that an $R$-module $G$
is torsion-free if $g s = 0 \; (g \in G, s \in S)$ only holds if $g = 0$. 
Moreover, $G$ is reduced (for~$S$) if
$\bigcap\limits_{s \in S} Gs = 0$. 
Throughout we suppose that $R$ is reduced and torsion-free (for~$S$).
The reader will observe that under these restrictions two kinds of 
realization theorems for $R$-algebras $A$ as endomorphism algebras of 
suitable modules $G$ are known. If we are lucky, then we find an $R$-module 
$G$ with 

$$\End_R G = A. \eqno(STRONG)$$

This first case we shall call a strong 
realization theorem. The constructed module is an $A$-module and 
multiplication by $a \in A$ is an $R$--endomorphism of $G$ because $R$ is 
commutative, hence $A \subseteq \End_RG$ where $\End_RG$ is the 
endomorphism ring of $G$ 
and the construction of $G$ shows how to get rid of the endomorphisms not 
in $A$.

The first deep result for a strong realization theorem is Corner's theorem 
\cite{Co1} mentioned at several places in this volume. Note that Corner at 
this time was interested in $R = \Z$ and $A$ torsion-free, reduced of 
cardinality $\aleph_0$ with special emphasis to rings $A$ of finite rank. 
Extensions and interesting applications of this result are due to Adalberto 
Orsatti \cite{O1,O2,O3} whom we want to honour by including
this paper into a volume of articles on module theory for his 60th birthday.
\smallskip

Corner's result was extended in a number of papers which we do not 
want to discuss in detail. The reader is asked to consult the `unified 
treatment' in
Corner, G\"obel \cite{CG} which extends known results and also summarizes 
the new developments in the early eighties. Paper \cite{CG} is based on new 
combinatorial techniques first used for $p$-groups in Shelah \cite{S3}, then 
refined in \cite{S4,S5,S6} to what is called after \cite{CG} `Shelah's 
Black Box'. We only mention some of the main contributions obtained 
\cite{Co1, Co2, Co3, CG, DG1,DG2, DG3, DG4, GM1, GM2, FG, S3, S4, S5, S6} 
and surveys in \cite{DG5, EM}.
\smallskip

Besides the case of strong realizations it will happen quite often that
mathematical interest leads to less lucky cases. We cannot expect a strong
realization theorem due to unavoidable endomorphisms. Classical
examples for this second kind of realization theorems 
are those for are abelian p-groups 
$G$ where it is known from early results that many 
small endomorphisms automatically 
belong to $\End G$, see Fuchs \cite{F1}. In 
order to derive a realization theorem for a `decent' ring (like the $p$-adic
integers $A = J_p$) we must replace (STRONG) by a weaker demand  
${\End G/\Small G \cong A}$,
which was investigated in Corner \cite{Co2} and for cardinals 
$\geq 2^{\aleph_0}$ in Shelah \cite{S3}, see Dugas, G\"obel \cite{DG1} for 
an extension. 

If we change the category from $p$-groups to some other class 
of modules, the ideal $\Small G$ must be replaced by some suitable  
ideal depending on that category. A useful definition of such an ideal 
should also reduce to well-known ideals like $\Small G$ or $\Fin G$ for 
well-studied categories.
 This idea was followed up in joint work with Dugas \cite{DG3}, in
\cite{DGG} and \cite{CG}, and lead to the ideal $\Ines G$ which (up to small
adjustments like purity) is the collection of all those 
$\sigma \in \End G$
which do extend to any 

$$\widehat{\sigma}: G' \arr G \; (G' \supseteq G)$$ 

where $G'$
belongs to the category under consideration. The desired weak realization
theorem is of the form

$$\End_R G = A \ltimes \Ines G  \eqno(WEAK)$$

where (WEAK) is a little stronger than 
$\End_R G/$ \Ines  $G \cong A$ and denotes a ring split extension. 

\medskip

The strong realization theorem is still a special case ($\Ines G = 0$) 
and can be obtained if $A$ is cotorsion-free, i.e. if 
$\Hom_R (\widehat{R}, A) = 0$ where $\widehat{R}$ denotes the completion of 
$R$ in the $S$-topology (generated by $Rs, s \in S)$. 
Cotorsion-free modules $G$ of cardinality $\lambda$ with (STRONG) 
have been constructed earlier, see \cite[p. 456]{CG}. They all satisfy

$$  |G| = \lambda > |A| \text{ with } \lambda^{\aleph_0} 
= \lambda.\eqno(CARD)$$ 

We note that such cardinals $\lambda$ are not cofinal to $\omega$ by 
K\"onig's lemma, hence cardinals like $\lambda = \aleph_\omega$ are 
excluded. The proof of the Black Box uses $\lambda^{\aleph_0} = \lambda$, so 
the restriction seems to be due to the Black Box. 
The same holds for weak realization theorems. 

It is the aim of this paper to study this drawback more closely. We want to 
deal with this in two more definite classical cases where either $\Ines G 
=0$ or $\Ines G = \Fin G $
where $\Fin G$ denotes the ideal of all $\sigma \in \End G$ with $\Im 
\sigma$ of finite rank.
The latter case comes up naturally for two classical 
categories, separable abelian groups and  $\aleph_0$-cotorsion-free modules.
Separable modules are pure submodules of products $\prod R$ of the 
ring $R$ and $\aleph_0$-cotorsion-free $G$ are
defined by the requirement that they are reduced and 
torsion-free such that every homomorphism from a complete module into $G$ 
has finite $p$-adic rank, see \cite{CG}.

\bigskip

If $A$ is $\aleph_0$-cotorsion-free, then we can find (e.g. in \cite 
[p. 470]{CG}) $\aleph_0$-cotorsion-free $R$-modules $G$ with 

$$(WEAK), \quad \Ines G = \Fin G \text{ and (CARD).} \eqno(FIN)$$

Similarly, if $A$ is $R$-free and countable, then we can find separable 
$R$-modules $G$ with (FIN); see \cite{DG4} and Corner, G\"obel \cite{CG1}.
If $A$ is uncountable, then we must add a technical condition discussed
in \cite{CG1,DG4}.

\medskip

We now come to our main concern, the problem whether we are able to avoid 
the cardinality
restriction on $\lambda$ caused (virtually) by the use of the Black Box
which is $\lambda = \lambda^{\aleph_0}$. 
The condition $\lambda^{\aleph_0} = \lambda$
is needed to complete an easy and transparent 
counting argument for predicting homomorphisms, see Appendix of 
\cite{CG}. 
Hence $\lambda^{\aleph_0} > \lambda$ requires at least changes of the Black 
Box. However, in Section 3 we will see that this obstacle is more basic 
and really not due to Shelah's Black Box but caused by the `natural 
algebraic' setting which is prepared for its use. Recall that the 
desired $R$-modules in all cases are sandwiched between a 
base module $B$ and its $S$-adic completion $\widehat{B}$, i.e.

$$B \subseteq G \subseteq \widehat{B}.$$

This initial step already removes the chance to work with $\lambda$ such that
${\cf (\lambda) = \omega}$ in case (WEAK) as follows from one of our main 
result:

\medskip

{\bf Corollary \ref{3.6}}
{\it Let $R$ be a principal ideal domain and $G$ be a torsion-free,
reduced $R$-module of cardinality $\lambda$ such that $\cf (\lambda) = \omega$.
Suppose $\mu^{|R|} < \lambda$ for all cardinals $\mu < \lambda$. If $G$ has
$\lambda$ pairwise distinct pure injective submodules, then $\End_R G/\Fin G$
has rank $\lambda^{\aleph_0}$.}

\bigskip

If we want to construct $\aleph_0$-cotorsion-free $R$-modules realizing an
$\aleph_0$-cotorsion-free (but not cotorsion-free) algebra $A$, then the
base module $B$ above is $\bigoplus\limits_{\lambda} A_R$ and $|A|< \lambda$.
Each copy of the $R$-module $A_R$ has a non-trivial cotorsion submodule and
$B$ as well as $G$, if of size 
$\lambda$, satisfies the requirements of Corollary
\ref{3.6} above. If $A \cong \End_R G/\Fin G$ then 
$|A| < \lambda$ contradicts the conclusion of
Corollary \ref{3.6}. Hence modules of cardinality $\lambda$ do not have the
desired endomorphism ring.
\bigskip

It is interesting to note that cofinality $\cf \lambda = \omega$ is used 
in the proof of (\ref{3.6}) to conclude

$$|\End_R G/\Fin G| = \lambda^{\aleph_0} > \lambda \text{ from }
|\End_R G| = \lambda^{\aleph_0} \text{ and } |\Fin G| \leq \lambda.$$

\bigskip

If $A$ is cotorsion-free, we have seen that the construction by the 
Black Box must be improved. 
In Section 4 we distinguish two cases (A) and (B) depending on the 
algebra $A$. In case (A) we assume that $A$ (as above) is 
cotorsion-free. A new combinatorial argument is introduced which is a 
mixture of the Black Box 
and an older combinatorial principle from \cite{S}, which 
was also used in G\"obel, May \cite{GM2} and 
named `Shelah elevator'. The 
Shelah elevator was used originally in \cite{S} for constructing arbitrary 
large indecomposable, torsion-free abelian groups. 
In \cite{GM2} it was used to export modules from smaller to larger 
cardinality. Here we first construct a 
fully rigid system of $R$-modules 
$G_X$ $(X \leq \mu)$ such that 
$$\Hom_R(G_X, G_Y) = A \mbox{ and } G_X \subseteq G_Y \mbox { if } 
X \subseteq Y \subseteq \mu$$ 
and 
$$\Hom_R (G_X, G_Y) = 0 \mbox{ if } X \not\subseteq Y.$$ 

If $\mu = \mu^{\aleph_0} < \lambda$ is some cardinal $ > |A|$, then the 
Black Box applies and we obtain a fully rigid system. 
In the second step we take a few members of this fully rigid 
system and put them into
the Shelah elevator and lift them up to $G$
with $|G| = \lambda$  as desired.

\bigskip

If $A$ is not cotorsion-free, we have to work harder to circumvent 
the dead end by (\ref{3.6}). We must avoid that 
$G$ has too many pure injective submodules. 
This is done in case (B).
Again, a basic idea is to carry information from a rigid system of 
$R$-modules
of smaller cardinal $\mu < \lambda$ up to $\lambda$. 
However, this time 
the Black Box is used to obtain an even 
stronger fully rigid system. 
Here a family of $R$-modules 
$\{G_u: u \in  \mu^{\leq \aleph_0}\}$ is called 
{\em essentially $A$-rigid over a 
directed subset $\mathfrak{U}$ of $\mu ^{\leq \aleph_0}$},
if the $G_{u}$'s are fully rigid as usually (see e.g. \cite{CG}):
$$\Hom_R(G_u, G_{u'}) = A \delta_{uu'}, \ltimes \Fin(G_u,G_{u'})$$
where 
$\delta_{uu'} = 1$ if $u \subseteq u'$ and 
$\delta_{uu'} = 0$ if $u \not\subseteq u'$.
Moreover $G_u \subseteq G_{u'}$ for $u \subseteq u' \in \mu^{\leq 
\aleph_0}$.\\
Hence $G_\mathfrak{U} = \bigcup\limits_{u \in \mathfrak{U}} G_u$ 
is a well defined $R$-module and
`rigidness' between any $G_u$ and $G_{\mathfrak{U}}$ is required as well; 
see Definition \ref{4.3}.\\
Inspection of the proofs in \cite[The torsion-free theory, pp. 464 - 
 465, $\Ines$ in other torsion-free theories (pp. 465 - 
470)]{CG} shows that the
existence of essentially rigid families can be replaced by these stronger
essentially rigid families, see (\ref{4.4}) and (\ref{4.5}). 
The main burden in the rest
of case (B) is to find a suitable directed system $\mathfrak{U}$ of size 
$\lambda$ to ensure that $G_{\mathfrak{U}}$ is of 
 size $\lambda$. Since we start from a 
family of modules of size $\mu$ given by the Black Box, $\lambda$ must be 
close enough to that $\mu$. If this is the case we derive a new realization 
theorem for algebras $A$ with particular emphasis on cardinals $\lambda$ 
with $\cf \lambda = \omega$.\\The main result is

\medskip

{\bf Theorem \ref{4.7}} {\it Let $A$ be an $R$-algebra, $\mu, \lambda$ be
cardinals such that $|A| \leq \mu = \mu^{\aleph_0} < \lambda \leq 2 ^\mu$. If $A$ is
$\aleph_0$-cotorsion-free or $A$ is countably free, 
respectively, then there exists an $\aleph_0$-cotorsion-free or a separable 
(reduced, torsion-free) $R$-module $G$ respectively of cardinality $|G| = 
\lambda$ with $\End_R G = A \oplus \Fin G$.}

\section{Basic definitions, examples and motivations}

Let $R$ be a principal ideal domain and $\chi = |R|^+$ be the successor of the
cardinality $|R|$ of $R$, which is fixed throughout Section 2 and 3.

\begin{definition} \label{2.1} We will say that an $R$-module $M$ of rank $\lambda \geq \chi$
has {\bf many pure injectives} if there are $\lambda$ pairwise distinct 
pure injective summands of $M$ (purely) generated by $< \chi$ elements.
\end{definition}

This definition may also be useful in the countable case as well, however
we are mainly interested in application close to Black Box proofs, hence
$\lambda \geq \aleph_1$.

\medskip
{\bf Examples}

Any module $M$ of rank $\lambda$ over a discrete 
valuation ring $R$ possesses a basic
submodule $B $, which is unique up to isomorphism; 
see Fuchs \cite{F1} or Eklof,
Mekler \cite[p.124]{EM}. 
Hence $B= \bigoplus_{i\in I} b_iR $ is a direct sum of pure cyclic submodules
$b_i R \; (i \in I)$. It is often the case that $M$ has many pure injectives:
\begin{enumerate}
\item If $R=J_p$ is the ring of 
$p$-adic integers, then $M$ is a direct sum of a divisible module
$D$ and a reduced submodule $M'$. If $D$ has rank $\lambda$, then $M$ 
has enough pure injectives. Otherwise we may 
assume that $D = 0$ and $M$ is a reduced $J_{p}$-module.
If $M$ is an abelian $p$-group, then a 
theorem of Kulikov applies, see Fuchs \cite[p. 146]{F1}. It shows 
when $|I| \ge \lambda$, then $M$ has many pure injectives.
If $M$ is torsion-free, then each summand of its basic module is
 pure injective. Hence $M$ has many 
injectives if, again, $|I| \geq \lambda$.
Also note that 
$B \subseteq M \subseteq \widehat{B}$ where $\widehat{B}$ is the $p$-adic 
completion of $B$; see Fuchs \cite{E1}.
\item The last remark relates to modules used in Black
 Box proofs for realizing rings as endomorphism rings; see Dugas, G\"obel
\cite{DG1, DG2, DG3}, Shelah \cite{S5, S6} or Corner, G\"obel \cite{CG}. In 
any case (mixed, torsion-free or torsion) - 
constructions begin with an $A$-submodule 
$B = \bigoplus\limits_{i \in \lambda}b_i A$  
of the final module $G$ with $\End G$ as required and
$$B \subseteq G \subseteq_* \widehat{B} \eqno{(*)}$$
where $\widehat{B}$ is 
the $S$-completion and $\subseteq_*$ denotes pure submodules. In all cases 
which are not cotorsion-free, the pure cyclic $A$-module $b_i A$ is not 
cotorsion-free and possesses a pure injective submodule $\neq 0$. Hence $G$ 
has many pure injectives if $G$ has rank $\lambda > |A|$. The final module
has size  $|G| = \lambda^{\aleph_0}$ 
which is $\lambda$ only if $\cf \lambda > \omega$. 
\end{enumerate}

We want to investigate what happens if we require 
${|G| = \lambda}$ and ${\cf \lambda = \omega}$.
Surely, many pure injectives may prevent the existence of realization 
theorems. Hence we consider this possibility first.

In  case $|G| = \lambda, \cf \lambda = \omega$ we note that 
(like in case $\cf \lambda >\omega)$,
the resulting module $G$ has many pure injectives. 
If however $|G| = \lambda^{\aleph_0}$ this is no harm. (In fact 
 it is not obvious from ($*$) 
and surprisingly not true as we shall show that $G$ (derived in the 
realization theorems) has many pure injectives.

\bigskip
If $\lambda$ has cofinality $\omega$, then we are 
bound to distinguish two cases. 
If the algebra $A$ is cotorsion-free, then we 
will derive new realization theorems for modules of size $\lambda$ cofinal to
$\omega$, which is similar to the known ones in  \cite{CG, S5,S6}.
If the module has many pure injectives, then we want to show that 
realization theorems (even modulo large ideals of inessential endomorphisms) 
do not exist. 
If the algebra $A$ is not cotorsion-free, in Section 4 we also find a 
way around  
to construct modules of size $\lambda$ with $\cf \lambda = \omega$ 
for certain cardinals having a specified endomorphism ring as before.

\section{Torsion-free $R$-modules - non existence of a realization theorem}

Recall that $R$ is a PID such that $R$ is reduced (and torsion-free) for 
some fixed multiplicatively closed, countable subset S.
Also recall that $N_*$ denotes the pure closure of 
$N \subset M$ if the $R$-module $M$ is torsion-free. Let $J_p$ denote the 
$p$-adic integers for some prime $p$, this is the $p$-adic completion of
$R$ provided $R$ is $p$-reduced.

We begin with a known result which appears in Dugas, G\"obel 
\cite{DG2}, see the proof in \cite{DG2} or in \cite{EM}. 

\bigskip

\begin{proposition} \label{3.1}
If $M$ is a reduced, torsion-free $R$-module and $N \subset M$
with $N \cong J_p$, then $N_* \cong J_p$ and $N_*$ is a summand of $M$.
\end{proposition}

\bigskip

\begin{corollary} \label{3.2}
Let $\lambda > |R|$ be some cardinal. If $M$ is torsion-free
reduced with pairwise distinct submodules $N_i \; (i \in \lambda)$ which are
pure--injective, then $M$ has $\lambda$ pure injectives which constitute a direct sum in $M$.
\end{corollary}

{\bf Proof}. Each $N_i$ is a $p$-adic module. We replace the given family by
an equipotent subfamily of $J_p$-modules for a fixed $p$. Similarly we may
assume that $N_i \cong J_p$. If we replace the new family by $(N_i)_* (i \leq \lambda)$,
each $(N_i)_*$ may coincide with finitely many $(N_j)_*$ by (\ref{3.1}). An
equipotent subfamily $J_p \cong N'_i \sqsubset M$ satisfies $\bigoplus\limits_{i \in \lambda} 
N'_i \subseteq M$. 
\hfill$\square$

\bigskip

The conclusion of the following Proposition \ref{3.3} follows from the
existence of a family of submodules similar to the one in \ref{3.2}. Under
these conditions it will be possible to find many endomorphisms.
These endomorphisms will destroy any hope
for a realization theorem, even modulo some ideal of 
inessential endomorphisms.
Moreover (\ref{3.3}) illustrates that (\ref{3.2}) must be strengthened in 
order to carry out (\ref{3.3}) and its consequences. Notice that (\ref{3.3}) 
is the main tool for proving the non-existence of a realization theorem.

\bigskip

\begin{proposition} \label{3.3}
Suppose $G = \bigcup\limits_{n \in \omega} G_n$ is the
union of a chain of pure submodules $G_n$ of cardinality 
$\lambda_n (n \in \omega)$
such that $\lambda_n \  (n \in \omega)$ is strictly increasing.
Let $\{N^n_i| i \in \lambda_{n+1}, n \in \omega \}$ be a family of 
pure injective modules such that $\bigoplus\limits_{i \in \lambda_{n+1}} 
N^n_i \subseteq G_{n+1}$ is direct and 
$G_n \oplus N^n_i \subseteq_* G_{n+1}$ 
is pure for any $n \in \omega$ and 
$i\in \lambda_{n+1}$. If $\eta \in \prod\limits_{n \in \omega} 
\lambda_{n+1}$
then there exists an $h_\eta \in \End_R G$ with 
$\Im h_\eta = \bigoplus\limits_{n \in \omega}N^n_{\eta(n)}$.
\end{proposition}

\begin{remark} {\rm Note that 
 $\lambda = |G| = \sup\limits_{n \in \omega} \lambda_n$ is cofinal to 
$\omega$.
The choice of $\Im h_\eta$ will ensure that $h_\eta$ is not
swallowed by $\Ines G$, the ideal of inessential endomorphisms of $G$. 
Obviously (\ref{3.3}) will imply the existence of 
$\lambda^{\aleph_0} > \lambda$
such endomorphisms and $R \cong \End G/\Ines G$ would be impossible for any
ring $R$ with $|R| \leq \lambda$; see (\ref{3.5}).}
\end{remark}

\medskip

{\bf Proof of (\ref{3.3})}. Let $h^n_i : G_n \oplus N^n_i \arr N^n_i$ be the
canonical projection. This projection extends to
$$h^n_i : G \arr N^n_i$$
because $G_n \oplus N^n_i$ is pure in $G$ and $N^n_i$ is pure injective. If
$\eta \in \prod
\limits_{n \in \omega} \lambda_{n+1}$, then put
$$h_\eta = \sum\limits_{n \in \omega} h^n_{\eta(n)}.$$
If $x \in G$, then there is $n \in \omega$ such that $x \in G_n$, hence
$h^m_{\eta (m)} (x) = 0$ for all $m \geq n$ and the sum
$h_\eta (x) = \sum\limits_{n \in \omega} h^n_{\eta (n)} (x)$ 
is finite and hence
well-defined in $G$. Clearly $h_\eta \in \End G$ and 
$\Im h_\eta = \bigoplus\limits_{n \in \omega} N^n_{\eta(n)}.$
\hfill$\square$

In view of (\ref{3.3}) we want to strengthen (\ref{3.2}).

\medskip

\begin{lemma} \label{3.4}
Let $\mu$ be a regular cardinal $> \chi$ and let $G$ be a torsion-free,
reduced $R$-module with the following properties.
\begin{enumerate}
\item There is a family $N_i \subseteq G$ $(i < \mu)$ of pure injective pairwise
distinct submodules purely generated by $< \chi$ elements.
\item Let $K \subseteq G$ with $|K| \leq \kappa < \mu$ for $s m e$ regular 
cardinal $\kappa$.
\end{enumerate}
Then we can find pure injective summands $0 \neq N'_i$ of $G$ and
$K' \subseteq G$ such that $K \subseteq K'$, $|K| \leq \kappa$ and
$K' \oplus \bigoplus\limits_{i < \mu} N'_i \subseteq_* G.$
\end{lemma}

{\bf Proof}.
By Corollary \ref{3.2} we replace the given family by a new family of pure 
injective summands $N_i\ne 0$ $(i < \mu)$
 such that $\bigoplus\limits_{i < \mu} N_i$ is a direct sum. 
Inductively we enumerate a subfamily of the
$N_i$'s and choose $K = K_0, K_i \subseteq K_{i+1}$ with $N_i \subseteq 
K_{i+1}$ and $N_i \oplus K_i$ which is a strictly increasing, continuous 
chain $K_i$ of submodules and elementary submodels of $G$ with respect to a 
language $L_\chi$ of cardinality $\chi$.

If $K_i$ is given, we want to find $N_i$ from the above family with 
$K_i \oplus N_i$.
Then we let $K_{i+1}$ be the elementary closure of $K_i \oplus N_i$ 
and proceed continuously. 
Recall that 
$J_p \cong N_j$ from (\ref{3.2}) and let $g_j : J_p \hookrightarrow G$ 
be the given isomorphism. 
There is some $g = g_j$ with $g(1) \notin K_i$ 
by cardinality. The algebraic reason for taking the elementary closure is 
that $G/K_i$ must be torsion-free, reduced. 
We want to show that $\Im g \cap K_i = 0$. 
Suppose $g(x) \in K_i$ for some $0 \neq x \in J_p$. If $x$ is not 
pure in $J_p$, then $x = p^k x'$ for some pure $x' \in J_p$. 
Hence $g (x) = p^k g (x') \in K_i$ and $K_i$ 
is pure in $M$ and torsion-free. We also have 
$g (x') \in K_i$ and hence may assume that $x$ is pure in $J_p$. 
There is a 
maximal $p$-power $p^k$ such that 

$p^k| g (1)$ modulo $K_i$ because $G/ K_i$ 
is torsion-free reduced and $0 \not\equiv g (1) + K_i$. 

We also find $n \in 
R$ such that $p^{k+1} | (n-x)$ in $J_p$, hence $p^{k+1} |(n g (1) - g (x))$ 
and $p^{k+1} | n g (1) \mod{K_i}$. We conclude $p|n$ and $p^{k+1} | (n-x)$ 
forces $p|x$ in $J_p$, contradicting purity. We have $N_i \oplus K_i$ for 
the above $N_j$ renamed as $N_i$. The above family $N_i, K_i$ $(i < \mu)$ is 
established.

Let $S = \{ \alpha < \mu: \cf (\alpha) \geq \chi\}$ which is stationary in
$\mu$. The following arguments do not use the specific structure of $N_i$. 
We only need that the $N_i$'s are the pure closure of $< \chi$ elements.
Also in the last paragraph we could have dropped the reference to 
(\ref{3.2}). If $\delta \in S$, then
$$K_\delta \oplus N_\delta \subseteq G$$
by the above family. The elementary submodel $K_\delta$ over $L_\chi$
ensures that $K_\delta$ allows an elementary embedding 
$h_\delta : N_\delta \arr K_\delta$.
Let 

$\Gamma_\delta$  be the set of equations $r | (x_i - c_r), (c_r \in 
K_\delta),(r \in R)$

where $x_i \in L \chi$ corresponds to some generator $a_i^\delta$ 
of $N_\delta$ (say $N_\delta$ is purely generated by a set 
$\{a^\delta_i : i \in I_\delta\}$
of size $< \chi$) such that $r | (a^\delta_i - c_r)$ in 
$M$ ($r \in R$). 
Then $\Gamma_\delta$
has $< \chi$ variables and $|R| < |R|^+ = \chi$ 
equations. The elementary embedding ensures some strong purity.
$$\text{If } r | (x_i - c_r) \in \Gamma_\delta \text{ then } 
r | (h_\delta (a^\delta_i) - c_r) \text{ in } G. \eqno{(*)}$$
Put $N'_\delta = \{x - h_\delta (x) : x \in N_\delta\} \subseteq N_\delta \oplus K_\delta \subseteq G$
and notice that 
$$ N_\delta \to N'_\delta (x \to x - h_\delta (x))$$ 
gives an isomorphism.
Clearly $N_\delta \oplus K_\delta = N'_\delta \oplus K_\delta$ by definition
of $N'_\delta$. Let $N''_\delta = (N'_\delta)_*$ and note that $N''_\delta$ is
pure injective as well by the above. Also note that $N''_\delta \cap K_\delta = 0$
and $N''_\delta \oplus K_\delta$ must be pure in $G$. The pure injective
module $N''_\delta$ is the first of our $\mu$ candidates needed in (\ref{3.4}). The
others show up by an easy combinatorial trick based on Fodor's Lemma, see
Jech \cite{Je}. Recall that $K_\delta \oplus N_\delta$ and 
$h_\delta : N_\delta \to K_\delta$
may be viewed as a regressive function and $S$ is stationary. Copies of
$N_\delta$ in $K_\delta$ can be enumerated by ordinals 
$< \delta$ as $\cf (\delta) \geq \chi$,
hence $h_\delta : N_\gamma \to N_\delta$ for $\gamma < \delta$. By Fodor's
Lemma there is a stationary subset $S_1 \subseteq S$ with 
$h_\delta (N_\delta) = N$
for some fixed $N$ and all $\delta \in S_1$.

Choose $K' \subseteq_* K_{i_{0}}$ for some $i_0 \in S_1$ (minimal) with
$N \cup K \subseteq K'$. 
Induction on $\delta$ with the last argument shows that
$$K' \oplus \bigoplus\limits_{\delta \in I} N''_\delta \subseteq_* G 
\mbox{ for } I = \{ \delta \in S_1, \delta \geq i_0\}.$$ 
Also note that $|I| = \mu$
and (\ref{3.4}) follows. \hfill$\square$

\bigskip

\begin{theorem} \label{3.5}
Let $R$ be a PID with $|R|^+ = \chi$ and let
$\lambda \geq \mu > \chi$ be cardinals with 
$\cf (\lambda) = \omega$ and
$\mu^{|R|} < \lambda$ for all $\mu < \lambda$.
If $G$ is a torsion-free, reduced $R$-module of cardinality $\lambda$ with a 
set of $\lambda$ pairwise distinct pure injective submodules, then we find 
$\lambda^{\aleph_0}$ endomorphisms $h_i (i \in \lambda^{\aleph_0)}$ with 
$\Im h_i$ pure and isomorphic to a direct sum of a countable infinite subset 
of some set of $\lambda$ pure injective submodules.
\end{theorem}

{\bf Proof}. Let $\lambda_n$ $(n \in \omega)$ be a strictly increasing
sequence of cardinals with $\sup\limits_{n \in \omega} \lambda _n = \lambda$.
Replacing $\lambda_n$ by its successor $\lambda^+_n$ if necessary, we may
assume that each $\lambda_n$ is a regular cardinal, moreover $\lambda_0 > \chi$.
We apply Lemma \ref{3.4} inductively to find a countable chain of pure submodules
$G_n \subseteq G$ such that $G = \bigcup\limits_{n \in \omega}G_n,$ $|G_n| = \lambda_n$
and such that there are pure injective modules $N^n_i$ $(i \in \lambda_{n+1})$ with
$$G_n \oplus \bigoplus\limits_{i \in \lambda_{n+1}} N^n_i \subseteq_* 
G_{n+1}.$$
Now we are in the position to apply Proposition \ref{3.3} and find 
endomorphisms 
$h_\eta \in \End G \;\;\; (\eta \in \prod\limits_{n \in 
\omega} \lambda_{n+1})$ such that 
$\Im h_\eta = \bigoplus\limits_{n \in \omega} N^n_{\eta(n)}$.
\hfill$\square$

\medskip

Realization theorems for
certain $R$-algebras $A$ provide $R$-modules $G$ with $A \cong \End G/J$ for
some suitable ideal $J = \Ines G$ depending on the nature of $A$ 
and modules $G$. If $G$ is torsion-free, then either $J = 0$ (in case $A$ 
is cotorsion-free) or $J = \Fin G$ is the ideal of those 
endomorphisms $\varphi$ of $G$ with $\Im \varphi$ of finite rank. More 
generally $J = \Ines G$ if $\Im \varphi$ is complete in the $S$-topology, 
see $\S$ 1 and 4.
\medskip

Our first application is an easy counting argument.

\medskip

\begin{corollary} \label{3.6}
Let $(R,G)$ be as in Theorem \ref{3.5}, then $\End G/\Fin G$ has
rank~ $\lambda^{\aleph_0}$.
\end{corollary}

\begin{remark} \label{3.8}
{\rm If the algebra $A$ is not cotorsion-free, then $A$ 
possesses a pure injective submodule $0 \neq N \subset A_R$ and 
any module $G$ in construction by the Black Box has a pure submodule 
$\bigoplus_\lambda A$, hence the hypothesis of (\ref{3.6})
holds, and $\End G/\Fin G \cong A$ 
is impossible.  This is 
in contrast to cardinals $\lambda$ with $\cf( \lambda) > \omega$, see 
\cite{CG} and \cite{FG}.}
\end{remark}

\medskip

{\bf Proof of \ref{3.6}}. Note that $|\End G| = \lambda^{\aleph_0}$ by
(\ref{3.5}) and $|\Fin G|=\lambda$, hence $|\End G/\Fin G| = \lambda^{\aleph_0}$
from $\lambda^{\aleph_0} > \lambda$.

\medskip

The next application is based on the observation that endomorphisms 
$h_i$ $(i\in \lambda^{\aleph_0})$ in Theorem \ref{3.5} are not complete: Each
$\Im h_i$ is pure and a countable direct sum of pure injectives. Hence
$h_i \notin \Ines G$ and a suitable choice of $h_i's$ ensures that the
following holds.

\bigskip

\begin{corollary} \label{3.7}
Let $(R,G)$ be as in Theorem \ref{3.5}, then $\End G/\Ines G$
has rank $\lambda^{\aleph_0}$ as well.
\end{corollary}

Sometimes the implication of Theorem \ref{3.5}
holds automatically, e.g. in case of certain classes of p-groups. In this
case (\ref{3.7}) follows by the given arguments. We leave it to the reader 
to check the details.

\bigskip

Remarks \ref{3.8} applies {\it mutatis mutandis} for Corollary \ref{3.7}. 
This 
might lead to the impressions that realization theorems (which so far have 
only been established for cardinals $\lambda$ with $\cf(\lambda) > \omega$ 
or if $R$ has `more than three primes') 
will always fail otherwise. Fortunately we will be able to extend 
the known results in Section 4.

\section{Realizing algebras}

Let $R$ be any fixed commutative ring, with a distinguished countable
multiplicatively closed subset  $S$ of non-zero-divisors as discussed in 
Section~1.  We will consider torsion-free, reduced $R$-algebras $A$ (for $S$).

\bigskip

In the first part (A) we concentrate on cotorsion-free $R$-algebras, so  we
require $\Hom_R (\widehat{R}, A) = 0$. Part (B) will be harder; we will  
deal with realization theorems of the (WEAK) form. 

As in Section 3  we choose a cardinal $\chi$ with $|A|^+ = \chi$.

\medskip

(A) In the cotorsion-free case we can follow an established road 
including only a little new work. However, we are mainly interested in 
cardinals $\lambda$ cofinal to $\omega$ and modules of this size, where
$\lambda > \chi$.

\medskip

\begin{theorem}. \label{4.1}
Let $A,R$ and $\chi$ be as above and suppose $\mu$ is a
cardinal with $\chi \leq \mu = \mu^{\aleph_0} \leq \lambda$. Then we can find
an $R$-module $G$ with $\End G = A$ and $|G| = \lambda$.
\end{theorem}
\medskip

We need a useful notion from  Corner \cite{Co3}, used in \cite{CG, 
GM1} and at many other places.
\medskip

\begin{definition} \label{4.2}
Let $A$ be an $R$-algebra. A family $\{G_X : X \subseteq I\}$
of $R$-modules $G_X$ is called fully $A$-rigid family over an indexing set
$I$ if for any subsets $X,Y, \subseteq I$ the following holds
$$\Hom_R (G_X, G_Y) = A \;\; \mbox{ and } \; G_X \subseteq G_Y 
\;\; \mbox{if} \;\; X \subseteq Y$$
$$\Hom_R (G_X, G_Y) = 0 \;\; \mbox{if} \;\; X \not\subseteq Y.$$
\end{definition}

\bigskip

{\bf Proof of (\ref{4.1}).} If $\lambda = \lambda^{\aleph_0}$, then 
the existence of a
fully $A$-rigid family over $\lambda$ follows from 
Corner, G\"obel \cite{CG} by
a proof based on Shelah's Black Box, see \cite{CG} and also \cite{S2}. In
particular, if $\lambda = \mu$, let $\{G_X : X \subseteq \mu \}$ 
be such a family.
Also note that $|G_X| = \mu$ follows from \cite{CG}. We choose a 
finite cotorsion-free rigid subfamily, taking a finite subset 
$I \subseteq  \mu$ with $|I| \geq 6$ and
$$\mathfrak{F} = \{G_X : X \subseteq I\}.$$

Note that $|I| = 4$ would suffice, see \cite{GM2}. This small family is the 
basic tool for applying a different combinatorial argument, 
the ``Shelah's elevator'',
see \cite{GM2} and also Shelah \cite{S1}. We will apply a version given in
Corner \cite{Co3} which can be used more directly
to obtain a cotorsion-free $R$-module
$G$ of cardinality $\lambda$ with $\End G = A$.
$\hfill{\square}$

\medskip

(B) In order to find realization theorems for algebras as 
endomorphisms algebras
$A$ of $R$-modules $G$ which have unavoidable inessential endomorphisms we
have to work harder for $\End G/ \Ines G \cong A$, where $\Ines G$ is the 
ideal of all
inessential endomorphisms of $G$. Since we are primarily interested in $G$'s
of cardinality $\lambda$ with 
$\cf (\lambda) = \omega$ we need different (new)
combinatorial techniques because the second combinatorial principle used in
(A) would break down. Nevertheless the new methods resembles ideas from this
method which originates from \cite{S}. While the proofs on this Shelah' 
elevator are based on 
a clever distribution of rigid pairs covering the forthcoming 
module, e.g. an indecomposable abelian group, the new method is no longer an 
elevator moving up from bottom to top (cardinals), see \cite{GM2}. It only 
connects certain levels, needs more fuel and runs on a more powerful rigid 
system (even more powerful then a fully rigid system), which we explain 
first. For clarity we restrict to modules $G$ with $\Ines G = \Fin G$, the 
unavoidable endomorphisms are those of finite rank. So we assume that the 
algebra $A$ is $\aleph_0$--cotorsion-free (e.g. $A = J_p)$, which 
automatically leads to $\Fin G$, see Corner, G\"obel \cite{CG}. Recall that 
an $R$-module $G$ is $\aleph_0$-cotorsion-free if $G$ is torsion-free, 
reduced and any cotorsion submodule has finite rank over $\widehat{R}$.

\begin{definition} \label{4.3} $\text{ }$
\begin{enumerate}
\item If $I$ is an indexing set of cardinality $\mu$, then
$J = P (I)^{\leq \aleph_0}$ denotes all subsets of cardinality $\leq \aleph_0$.
Obviously $J$ is partially ordered by inclusion and we will abuse notation
and write $\{i\} = i$ $(i \in I)$ for singletons.

\item Let $A$ be an $R$-algebra and $\mathfrak{U}$ be a directed subset of 
$J$. A family of $R$-modules $\{G_u : u \in J\}$ will be called an 
essentially $A$-rigid family for $\mathfrak{U}$ (over $\mu^{\leq \aleph_0})$ 
if the following holds.

\begin{enumerate}

\item If $u = \{u_i : i \leq n\}$ and the $u_i$'s in $J$ are pairwise 
disjoint, then $G_u = \bigoplus\limits_{i \leq n} G_{u_i}$

\item $\{G_u : u \in J\}$ is directed, if $u \subseteq u'$ then 
$G_u \subseteq _* G_{u'}.$
Let $G_{\mathfrak{U}} = \bigcup\limits_{u \in \mathfrak{U}} G_u$.

\item If $u \subseteq u'$, then $G_u \subseteq G_{u'}$ and 
$\Hom_R (G_u, G_{u'}) = A \oplus \Fin (G_u, G_u').$

\item If $u \in J$, then $\Hom_R (G_u, G_{\mathfrak{U}}) = A \oplus \Fin 
(G_u, G_{\mathfrak{U}}).$

\item If $u \not\subseteq u'$, then $\Hom_R (G_u, G_{u'}) = \Fin (G_u, G_{u'}).$

\item $|G_u| = \mu^{\aleph_0}$ for all $u \in J.$

\end{enumerate}
\end{enumerate}
\end{definition}

\medskip

An easy modification of the proof of 
the Main Theorem in \cite{CG} shows that we can
strengthen this result to get 

\bigskip

\begin{proposition} \label{4.4}
Let $|I| = \mu$ be a cardinal with $\mu^{\aleph_0} = \mu$
and $J, \mathfrak{U}$ be as in (\ref{4.3}). Let $A$ be an 
$\aleph_0$-cotorsion-free $R$-algebra with $|A| \leq \mu$. \\
Then we can find an essentially $A$-rigid family 
$\{G_u : u \in J\}$ of $\aleph_0$-cotorsion-free $R$-modules for 
$\mathfrak{U}$.
\end{proposition}

\medskip

{\bf Proof.} By inspection of \cite{CG}.

\bigskip

A similar result holds for separable modules. Separable $R$-modules are
submodules of products $\prod R$.

\medskip

\begin{proposition} \label{4.5}
Let $|I| = \mu$ be a cardinal with
$\mu^{\aleph_0} = \mu$ and $J, \mathfrak{U}$ be as in (\ref{4.3}).
Let $A$ be a free $R$-algebra 
which is countably generated or satisfies a technical
(`nasty') condition discussed in \cite{CG1, DG3} with $|A| < \mu$.\\
Then there exists an essentially rigid family $\{G_u : u \in J\}$ of separable
$R$-modules for $\mathfrak{U}$.
\end{proposition}
\medskip

{\bf Proof.} See \cite{DG4} and \cite{CG1} for $A$ countably generated.

The modules in (\ref{4.4}) and (\ref{4.5}) give rise to the desired modules 
$G_{\mathfrak{U}}$ for a suitable directed system $\mathfrak{U}$. The 
relevant properties of $\mathfrak{U}$ are derived in our next

\medskip

\begin{proposition} \label{4.6}
Let $W \subset [ \lambda]^{\leq \aleph_0}$ with
$|W| \leq \lambda$ and $\lambda^{\aleph_0} > \lambda > 
\kappa = \cf (\kappa) > \aleph_0$.
Then we can find a directed subset 
${\cal U} \subseteq [\lambda]^{\leq \aleph_0}$
and a coding function $\|\;\|: \lambda \ra {\cal U}$ 
with the following properties

\begin{enumerate}
\item $\|\;\|$ is a bijection
\item $W \cup [\lambda]^{\leq \aleph_0} \subseteq {\cal U}$
\item If $\alpha: \kappa \ra \lambda$, then there exists $u \in {\cal U}$ 
with
$|\{ i \in \kappa: \| \alpha (i) \| \subset u\} | \geq \aleph_0$.
\end{enumerate}
\end{proposition}

\medskip

{\bf Proof.} Write $\lambda = \overset{\cdot}{\bigcup\limits_{\zeta < 
\theta}} A_\zeta$ for a decomposition of $\lambda$ into $\theta$ subsets 
$A_\zeta$ of size $\lambda$ for some regular cardinal $\theta$ with $\lambda 
> \theta > \kappa$. \\
Let $V_\zeta = \bigcup\limits_{\beta < \zeta} A_\beta$ 
and $U_0 = W \cup [\lambda]^{< \aleph_0},\|\;\|_0 = \emptyset$.\\
We want to define inductively $U_\zeta \subseteq [\lambda]
^{\leq \aleph_0}, \|\;\|_\zeta : V_\zeta \ra U_\zeta$ for each $\zeta \in \theta$ as
ascending, continuous chains such that ${\cal U} = \bigcup\limits_{\zeta 
\in \theta} U_\zeta$ and $\|\;\| = \bigcup\limits_{\zeta \in \theta} 
\|\;\|_\zeta$. From $\dom \|\;\|_\zeta = V_\zeta$
follows  $\dom \|\;\| = \lambda$ immediately. Note that $|U_0| = \lambda$, 
and a bijective map $\|\;\|_1:A_0 \ra U_0$ can be defined. 
Suppose $\xi \in 
\theta$ and 
$\|\;\|_\zeta : V_\zeta \ra U_\zeta$ is defined for all $\zeta < \xi$. 
If $\xi$ is a limit ordinal we take unions. Suppose $\xi = \zeta + 1$, then 
we must define $\|\;\|_\xi$ and $U_\xi$ such that (a) and (c) hold for $U = 
U_\xi$ and those $\alpha$ with $\Im \alpha \subseteq V_\zeta$. Note that 
$\alpha$ is regular, so any $\alpha$ gives rise to some such $\zeta$. Hence 
(a) and (c) will also hold for $U$. Condition (a) requires only that 
$\|\;\|_\zeta$ is extended to $\|\;\|_\xi$ as a bijection $\|\;\|_\xi : 
V_\xi \ra U_\xi$. It remains to define any bijection $\|\;\| : A_\zeta 
\ra U_\xi \setminus U_\zeta$ taking care of (c).

\medskip

\relax From $|A_\zeta| = \lambda$ we need $|U_\xi \setminus U_\zeta| = \lambda$,
$U_\xi \supset U_\zeta$ and define $U_\xi$ in three steps.

\medskip

First we take the ideal $U'_\xi$ generated by $U_\xi$, that is
$$U'_\zeta = \{ u \subset \bigcup E : E \subset U_\zeta, E \text{ finite}\}.$$
Then adjoin the $\|\;\|_\zeta$-closure set of $U_\zeta$, the set  
$\bar{U}_\zeta$ of all elements
$$u' = \bigcup \{ \| i \|_\zeta : i \in u \cap V_\zeta \}$$
where $u \in U_\zeta$. 
Note that $U''_\zeta = U'_\zeta \cup \bar{U}_\zeta$ has
size $\lambda$ while $\lambda^{\aleph_0} > \lambda$,
hence $|[\lambda]^{\leq \aleph_0} \setminus
U''_\zeta | = \lambda^{\aleph_0}$ and we also find a set 
$U'''_\zeta \subseteq [\lambda]^
{\leq \aleph_0} \setminus U''_\zeta$ of cardinal $\lambda$. Now we define 
$U_\xi = U'''_\zeta \cup U''_\zeta$ and $U$ is constructed.

\medskip

The first step in the construction of $U_\xi$ ensures that $U$ 
is directed in
$$([\lambda]^{\leq \aleph_0}, \subseteq)$$ 
and the second step ensures (c).
If $\alpha : \kappa \ra \lambda$ then $\Im \alpha \subseteq V_\zeta$ for some
$\zeta < \theta$.
Next we consider $\|\Im \alpha \| = \{ \| \alpha (i) \| : i \in \kappa \} \subseteq U$.
The closure properties provide $u \in U$ with $\| \alpha (i) \| \subseteq U$
for infinitely many $i \in \kappa$.

\bigskip

\begin{theorem} \label{4.7}
Let $A$ be an $R$-algebra and $\mu, \lambda$
cardinals such that $|A| \leq \mu = \mu^{\aleph_0} < \lambda \leq 2^\mu$.

\begin{enumerate}
\item If $A_R$ is $\aleph_0$-cotorsion-free, then there exists an 
$\aleph_0$-cotorsion-free $R$-module $G$ of cardinality $\lambda$ with 
$\End_R G = A \oplus \Fin G$.
\item If $A_R$ is a free $R$-module which is
either countably generated or satisfies the `nasty' condition from
\cite{DG4}, then there exists a separable (torsion-free)
$R$-module $G$ of cardinality $\lambda$ with $\End_R G = A \oplus \Fin G$.
\end{enumerate}
\end{theorem}

\medskip

\begin{remark}
{\rm Theorem \ref{4.7} is new for cardinals $\lambda$ cofinal to $\omega$.}
\end{remark}

\medskip

{\bf Proof.} Let $\{X_i : i \in \lambda\}$ be an anti-chain of the power set 
$P(\mu)$ and choose $U$ from Proposition \ref{4.6}. Then we build the new 
partially ordered set
$$\mathfrak{U} = \{\{X_i : i \in u\} : u \in {\cal U}\}.$$
By Proposition \ref{4.4} or
Proposition \ref{4.3}, respectively there is an essentially rigid family
$\{G_u : u \in [\mu]^{\leq \aleph_0}\}$ for $\mathfrak{U}$ of 
$\aleph_0$-cotorsion-free $R$-modules or separable $R$-modules respectively. 
In particular $G = G_\mathfrak{U} =\bigcup\limits_{u \in U} G_u$.

If $\sigma \in \End_R G$ and $u \in \mathfrak{U}$ then there exists $a_u \in 
A$ such that
$$(\sigma \upharpoonright G_u) - a_u \in \Fin_u G$$
from Definition \ref{4.3}.
Here we identify $a_u \in A$ with scalar multiplication $a_u$ on $G$ (or any
submodule). Moreover $\Fin_u G$ 
denotes (partial) homomorphisms from $G_u$ into
$G$ of finite rank. It follows 
immediately that $a_u$ does not depend on $u$ since
$G$ has infinite rank, is torsion-free and $\mathfrak{U}$ is directed:\\
$(\sigma \upharpoonright G_u) - a_u, (\sigma \upharpoonright G_{u'}) - a_{u'}
\subseteq (\sigma \upharpoonright G_v) -a_v$ for some $v \in \mathfrak{U}$ 
implies $a = a_u = a_v = a_{u'}$, hence

$$(\sigma \upharpoonright G_u) - a \in \Fin_u G \text{ for all } u \in 
\mathfrak{U}. \eqno(*)$$

We also claim that $\sigma - a \in \Fin G$. Otherwise $\Im (\sigma-a)$ has
infinite rank, there are independent elements $y_i \in \Im (\sigma-a)$, $(i \in \omega)$.
If $\sigma' = \sigma - a$, then we find $x_i \in G_\mathfrak{U}$ 
with $x_i \sigma' = y_i$
$(i \in \omega)$. We may assume $x_i, y_i \in G_{u_{i}}$ with 
$u_i \in \mathfrak{U}$ $(i \in \omega)$.\\
By Proposition \ref{4.6} (a) there are $\alpha_i \in \lambda$ such
that $\| \alpha_i \| = u_i (i \in \omega)$. Moreover, by Proposition \ref{4.6} (c)
we can find $u \in U$ such that $\{ i \in \omega : \| \alpha_i \| \subseteq u\}$
is infinite. Hence $x_i, y_i \in G_{u_{i}} \subseteq G_u$ for infinitely 
many $i \in \omega$ from (\ref{4.3}) (ii). However these 
$y_i$ belong to 
$\Im((\sigma - a) \upharpoonright G_u)$ and are independent. The mapping 
$(\sigma - a) \upharpoonright G_u$ has infinite rank and 
$(\sigma - a )\upharpoonright G_u \not\in \Fin_u G$ contradicts $(*)$. Hence 
(a) and (b) follows.
\hfill$\square$

\bigskip

Finally we want to state a result for the kind of cardinal $\lambda$ not
covered by (\ref{4.7}).

\begin{theorem} \label{4.8}. Let $A$ be an $R$-algebra and $\lambda > |A|$ be
a strong limit cardinal of cofinality $\omega$. Then (a) and (b) from (\ref{4.7})
hold.
\end{theorem}
\medskip

{\bf Sketch of a proof.} We can write $\lambda = \sup\limits_{n \in \omega} \lambda_n$
such that $2^{\lambda}_{n} < \lambda_{n+1}$ and $\lambda_n$ is a successor
cardinal. Then we can use essentially rigid families for each cardinality
$\lambda_n$ using $\Diamond$-arguments, which give rise to the desired $R$-modules of
size $\lambda$. \hfill$\square$

\end{document}